\documentclass[reqno]{amsart}
\usepackage{amssymb}
\usepackage{amsfonts}
\usepackage{mathrsfs}
\usepackage{graphicx}
\usepackage{booktabs}
\usepackage{subfigure}
\usepackage{arydshln}
\usepackage{cite}
\usepackage{color}

\newtheorem{theorem}{Theorem}
\theoremstyle{plain}

\newtheorem{claim}{Claim}

\newtheorem{corollary}{Corollary}

\newtheorem{definition}{Definition}

\newtheorem{lemma}{Lemma}

\newtheorem{proposition}{Proposition}
\newtheorem{remark}{Remark}

\numberwithin{equation}{section}

\makeatletter
\def\subsection{\@startsection{subsection}{2}%
  \z@{.5\linespacing\@plus.7\linespacing}
  {1pt}%
  {\normalfont\bfseries}}
\makeatother

\begin{document}
\title[Scaling Limits and Fractal Dimensions of $\omega $-Substitution Graph Systems]%
{Scaling Limits and Fractal Dimensions of $\omega $-Substitution Graph Systems}
\author{Qingcheng Zeng}
\address{School of Mathematical Sciences, Beihang University, Beijing
100083, P. R. China}
\email{qczeng@buaa.edu.cn}
\author{Yumei Xue}
\address{School of Mathematical Sciences, Beihang University, Beijing
100083, P. R. China}
\email{yxue@buaa.edu.cn}
\author{Cheng Zeng}
\address{School of Mathematics and Information Science, Shandong Technology
and Business University, Yantai, Shandong Province 264003, P. R. China}
\email{czeng@sdtbu.edu.cn}
\thanks{Yumei Xue is the corresponding author.}

\begin{abstract}
This paper introduces a class of non-homogeneous self-similar constructions, called $\omega $-substitution graph system, in which the interface vertices are
determined by an infinite sequence $\omega .$ The associated graph sequences naturally form compact metric measure spaces, and under suitable assumptions on $\omega 
$, we establish their convergence in the Gromov--Hausdorff--Prokhorov sense. We determine the
Hausdorff and box-counting dimensions of the limiting space. We further prove that the discrete
Hausdorff dimension of the graph sequence coincides with the Hausdorff
dimension of its scaling limit. As a direct
application of the scaling limit, asymptotic formulas are derived for a class
of generalized Wiener indices, which include the classical Wiener index and
Kirchhoff index as special cases. Finally, closed-form expressions for these
indices are also obtained.
\end{abstract}

\maketitle

\vspace{-1.8em}
\noindent\hspace*{3pc}%
\begin{minipage}{\dimexpr\textwidth-6pc\relax}
\end{minipage}
\par\medskip

\section{Introduction}

Graph theory and fractal geometry have both undergone extensive development
since the twentieth century, and their interaction has become a central
theme in the study of complex networks. Motivated by empirical findings on
small-world networks \cite{Watts}, scale-free networks \cite{Bara} and
fractal networks \cite{Song}, it is now widely recognized that many
real-world systems manifest underlying structural regularities.

However, a systematic mathematical understanding of complex networks is
still incomplete, as most existing studies focus on empirical data or finite
graph realizations rather than infinite limiting objects. This naturally
motivates the search for limiting objects that describe the asymptotic
geometry of complex networks.

One natural approach to constructing such limiting objects is through
scaling limits of graph sequences. By equipping each graph with appropriate
normalized metric, one may regard the resulting sequence as a family of
compact metric spaces. The Gromov--Hausdorff topology provides a standard
framework for investigating the convergence of such spaces and has found numerous
applications in geometric analysis and random geometry \cite%
{Gromov,Burago,Aldous}. To accommodate richer structures arising in various
applications, several extensions have subsequently been developed, including
the Gromov--Hausdorff--Prokhorov topology for compact measure spaces \cite%
{Miermont,Abraham} and the Gromov--Hausdorff--Prokhorov--Uniform topology
for spatial stochastic processes \cite{Gwynne}. These topologies have become
fundamental tools in the study of scaling limits of random graphs \cite%
{Bhamidi} and stochastic processes \cite{Croydon}.

Substitution provides a natural mechanism for generating recursively defined
graph sequences with self-similar structures. The concept of graph
substitution was systematically formalized by Previte \cite{Previte}, who
established a general framework for replacing vertices or subgraphs with
finite graphs and investigated the dynamics of substitution together with
fractal dimensions of limiting spaces. Since then, substitution-based graph
constructions have attracted increasing attention in both graph theory and
fractal geometry.

Inspired by substitution and iterated function systems, Xi et al. \cite{Xi1}
introduced substitution networks generated by directed edge-substitution
rules, and established their fractality and scale-freeness. Subsequently,
Li, Ye and Xi \cite{Li-T} investigated the Markov spectra of substitution
networks. This framework was later extended to multi-colored structures by
Li, Yao and Wan \cite{Li1}, while the average distance behavior was studied
by Ye and Xi \cite{Ye}. Zeng et al. further studied discrete fractal
dimensions for substitution networks in two cases: when each seed graph is a
cycle \cite{Zeng1}, and when each seed graph is a tree \cite{Zeng2}.
Moreover, Li et al. \cite{Li2} investigated scale-freeness in the random
setting of colored substitution networks. More recently, such
substitution-based graph constructions have been used to generate a wide
variety of fractal spaces \cite{Anttila}. Moreover, fractal dimensions
of the associated Gromov--Hausdorff limiting spaces have recently been
investigated \cite{Li3}, while whether these dimensions can be recovered
from the corresponding graph sequences remains largely unexplored.

Motivated by the above substitution-based models, we introduce a class of $%
\boldsymbol{\omega }$\textbf{-substitution graph systems} in which the
interface vertices are determined by an infinite sequence $\omega $ rather
than fixed interface vertices, allowing the resulting graph sequences to
exhibit non-homogeneous structures.

In this paper, we study these graph sequences within the framework of rooted
metric measure spaces, and establish their convergence in the
Gromov--Hausdorff--Prokhorov sense under suitable assumptions on $\omega .$
Our construction is based on the metric quotient of the symbolic space,
providing an alternative approach to the limiting space. Its metric depends
not only on the longest common prefix but also on the suffix relative to the
distinguished interface $\omega $. We further prove a stability result for
fractal dimensions, showing consistency between the discrete Hausdorff of the
graph sequence and the Hausdorff dimension of the limiting space.

Distance-based graph invariants have also been investigated in self-similar
and substitution-based network models; see, for example, \cite%
{Lu,Wang,Wu,Jin}. The scaling limit further yields the leading-order
asymptotics of a class of generalized Wiener indices, encompassing both the
Wiener index and the Kirchhoff index as special cases, for which closed-form
expressions are also obtained.

\bigskip

\section{Main results}

\subsection{$\protect\omega $-substitution graph systems}

Let $\Sigma =\{1,2,\ldots ,N\}$ be a finite alphabet with $N\geq 2$ , and let $\mathbb{N}=\{1,2,\ldots\}$. An
element $\alpha \in \Sigma ^{\mathbb{N}}$ is called an \emph{infinite word}
which is a mapping from $\mathbb{N}$ to $\Sigma $. It is typically
represented as an infinite sequence $\alpha _{1}\alpha _{2}\cdots ,$ where
$\alpha _{i}=\alpha (i)$ for $i\in \mathbb{N}$.  On the other hand, an element $u$ in $\Sigma
^{\ast }:=\cup _{k=0}^{\infty }\Sigma ^{k}$ is called a \emph{finite word}
and represented as a finite sequence $u_{1}u_{2}\cdots u_{k}$ with $u_{i}\in
\Sigma ,$ where $k$ is such that $u\in \Sigma ^{k},$ which is called the 
\emph{length} of $u$ and is denoted by $\left\vert u\right\vert .$

For any $u\in \Sigma ^{\ast }$ and $\alpha \in \Sigma ^{\mathbb{N}},$ the 
\emph{concatenation} $u\alpha $ is the infinite word defined by $%
u_{1}u_{2}\cdots u_{k}\alpha _{1}\alpha _{2}\cdots ,$ where $k=\left\vert
u\right\vert .$ In this case, $u$ is called a \emph{prefix} of the infinite
word $u\alpha .$ Similarly, the concatenation $uv$ of $u$ and a finite word $%
v=v_{1}v_{2}\cdots v_{l}$ is defined as the finite word $u_{1}u_{2}\cdots
u_{k}v_{1}v_{2}\cdots v_{l};$ in this case, $u$ is a prefix of $uv.$ For $%
x\in \Sigma ^{\ast }\cup \Sigma ^{\mathbb{N}},$ we let $x|_{k}$ denote its
prefix of length $k$ (with $k\leq \left\vert x\right\vert $ if $x$ is
finite). In particular, for any non-empty finite word $u\in \Sigma ^{k},$ we
denote by $u^{-}$ its prefix of length $k-1$, where $u^{-}$ is the empty
word $\epsilon $ if $k=1$. For $u\in \Sigma ^{\ast },$ the cylinder set $[u]$
determined by $u$ is defined by%
\begin{equation*}
\lbrack u]=\{\alpha \in \Sigma ^{\mathbb{N}}:u\text{ is a prefix of }\alpha
\}.
\end{equation*}

For any two words $x,y\in \Sigma ^{\ast }\cup \Sigma ^{\mathbb{N}},$ we
denote by $x\wedge y$ their \emph{longest common prefix,} and $\left\vert
x\wedge y\right\vert $ represents its length. For any letter $a\in \Sigma ,$
we let $\overline{a}$ denote the constant infinite word $aaa\cdots ,$ in
which every letter is $a.$

\begin{definition}
Let $G$ be a connected simple graph with vertex set $\Sigma .$ Given an
infinite word $\omega \in \Sigma ^{\mathbb{N}}$, the sequence of $\omega $%
-substitution graphs $\{G_{n}^{\omega }\}_{n\geq 1}$ is defined recursively
as follows. Let $G_{1}^{\omega }=G$ and $V(G_{n}^{\omega })=\Sigma ^{n}$ for
all $n.$ For any $n>1$ and any two distinct vertices $u,v\in V_{n},$ $%
(u,v)\in E(G_{n}^{\omega })$ if either
\end{definition}

\begin{enumerate}
\item[(i)] $u^{-}=v^{-}$ and $(u_{n},v_{n})\in E(G),$ or

\item[(ii)] $(u^{-},v^{-})\in E(G_{n-1}^{\omega })$ and $u_{n}=v_{n}=\omega
_{n-\left\vert u\wedge v\right\vert -1}.$
\end{enumerate}

\begin{figure}[tbph]
\includegraphics[width=0.6\textwidth]{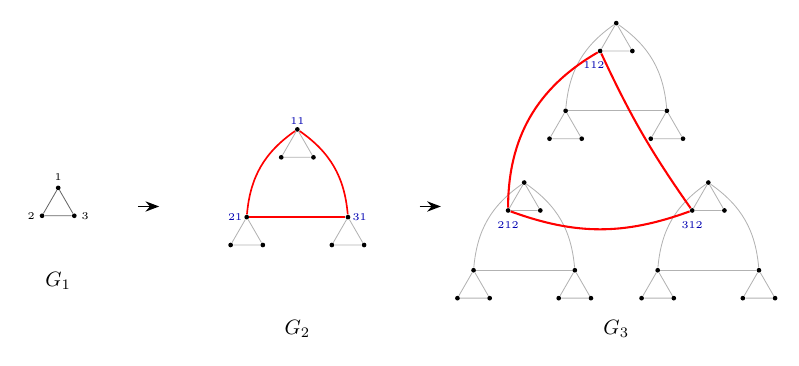} 
\caption{Evolution of the $\omega$-substitution graph system starting from $K_3$ with $\omega = 12\cdots $.}
\label{fig1}
\end{figure}

For notational simplicity, we shall write $G_{n}$ instead of $G_{n}^{\omega
} $ whenever no confusion can arise. An example of the $\omega$-substitution process for the case where $G$ is a 3-cycle $K_3$ and $\omega = 12\cdots $ is shown in Figure \ref{fig1}.

Let $d_{n}$ be a metric on $G_{n}$ for every $n\geq 1.$ The sequence $\{(G_{n},d_{n})\}$ is called an $\omega $\emph{-substitution graph
system}.

\begin{definition}
An $\omega $-substitution graph system $\{G_{n},d_{n}\}$ is called \emph{%
metric-preserving} if, for every $n\geq 1,$
\end{definition}

\begin{enumerate}
\item[(i)] $d_{n}(iu,iv)=d_{n-1}(u,v)$, $i\in \Sigma $ and $u,v\in \Sigma
^{n-1}$;

\item[(ii)] $d_{n}$ satisfies the cut-point additivity property, that is, if 
$z$ is a cut vertex of $G_{n}$ separating $u$ and $v,$ then 
\begin{equation*}
d_{n}(u,v)=d_{n}(u,z)+d_{n}(z,v);
\end{equation*}

\item[(iii)] $d_{n}(i(\omega |_{n-1}),j(\omega |_{n-1}))=d_{1}(i,j),$ $%
i,j\in \Sigma .$
\end{enumerate}

\begin{remark}
The metric-preserving property is satisfied by several natural metrics on an $%
\omega $-substitution graph system. In particular, if $d_{n}$ is the
geodesic distance or the effective resistance metric on $G_{n},$ then $%
\{(G_{n},d_{n})\}$ is metric-preserving.
\end{remark}

\subsection{Main theorems}

If $\omega $ is a constant sequence, then the substitution process
degenerates into a trivial recursive pattern. To avoid such trivial cases,
throughout this paper we assume that%
\begin{equation}
\omega _{k}\neq \omega _{1},\text{ }\forall k\geq 2.  \label{assump-omega}
\end{equation}%
In particular, $\omega $ is not periodic, i.e.,  $\sigma ^{t}\omega \neq \omega \text{ for all }t\in \mathbb{N}$, where $\sigma $ is the shift.

This condition guarantees a nontrivial recursive geometry, in particular,
(see Proposition \ref{prop-diam})%
\begin{equation*}
\mathrm{diam}(G_{n},d_{n})\asymp 2^{n}.
\end{equation*}%
Here notation $a_{\lambda }\asymp b_{\lambda }$ means that $a_{\lambda }$ is
comparable to $b_{\lambda }$ in a given index set $\Lambda $, i.e., there
exists a constant $C\geq 1$ such that $C^{-1}b_{\lambda }\leq a_{\lambda
}\leq Cb_{\lambda }$ for all $\lambda \in \Lambda .$

Motivated by this scaling behavior, for each $n\geq 1,$ we define the
normalized metric%
\begin{equation*}
\widetilde{d}_{n}:=2^{-n}d_{n}.
\end{equation*}%
For each $n\geq 1,$ we equip the finite symbolic space $\Sigma ^{n}$ with
the uniform discrete probability measure $\mu _{n},$ given by $\mu
_{n}(v)=N^{-n}$ for all $v\in \Sigma ^{n}.$

Let $\mathcal{X}_{n}^{\omega }=(\Sigma ^{n},\widetilde{d}_{n},\overline{1}%
|_{n},\mu _{n})$. The precise framework of rooted weighted metric space and
the associated Gromov--Hausdorff--Prokhorov topology will be recalled in
Subsection \ref{subsection-GHP}. Our first result states the convergence of $\{%
\mathcal{X}_{n}^{\omega }\}$ in this topology (also see Theorem \ref{th-GHP-1}).

\begin{theorem}
Let $\{(G_{n},d_{n})\}_{n\geq 1}$ be a metric-preserving $\omega $%
-substitution graph system. Then the associated sequence of compact rooted weighted metric space $\{\mathcal{X}_{n}^{\omega }\}$ converges in the
Gromov-Hausdorff-Prokhorov topology.
\end{theorem}

\begin{remark}
The above result guarantees the existence of a limiting compact rooted weighted metric space $\mathcal{X}_{\omega }=(X_{\omega },d_{\omega
},\emptyset _{\omega },\mu _{\omega })$ whose construction is given in
Section 4. The metric $d_{\omega }$ incorporates additional dependence on
the distinguished interface $\omega $, so that the resulting space differs
from standard self-similar symbolic spaces.
\end{remark}

The scaling limits obtained in the previous result provide a family of
compact metric spaces indexed by $\omega .$ We now investigate their fractal
geometry through two of the most fundamental notions of dimension, namely
the Hausdorff dimension and box dimension. Our first result shows that, although the
limiting spaces generally depend on the choice of $\omega ,$ these
dimensions are universal across the family and admit a simple explicit
formula determined solely by the substitution rule.

\begin{theorem}
\label{th-dim}Let $\mathcal{G}:=\{(G_{n},d_{n})\}_{n\geq 1}$ be an $\omega $%
-substitution graph system, where $d_{n}$ is the geodesic distance on $G_{n}$%
. Let $(X_{\omega },d_{\omega })$ be the limiting metric space associated
with $\mathcal{G}$. Then 
\begin{equation*}
\dim _{\mathrm{H}}(X_{\omega })=\dim _{\mathrm{B}}(X_{\omega })=\frac{\log N%
}{\log 2}.
\end{equation*}
\end{theorem}

The explicit dimension formula obtained above naturally raises the question
of whether this fractal dimension is already encoded in the approximating
graph sequence. We next show that the discrete Hausdorff dimension of the $%
\omega $-substitution graphs coincides with the Hausdorff dimension of the
associated scaling limits. The precise definition of the discrete Hausdorff
dimension will be given in Section 5.

\begin{theorem}
\label{th-dim-stab}Under the same assumptions of Theorem \ref{th-dim}, 
\begin{equation*}
\dim _{\mathrm{DH}}(\mathcal{G})=\dim _{\mathrm{H}}(X_{\omega }).
\end{equation*}
\end{theorem}

Having established the dimension result, we now turn to a global geometric
quantity describing the average metric structure of the limit space. Let 
\begin{equation*}
\mathcal{D}_{\omega }:= \int_{(x,y)\in X_{\omega }\times X_{\omega }}d_{\omega
}(x,y)d\mu _{\omega }(x)d\mu _{\omega }(y)
\end{equation*}%
be the \emph{average }$d_{\omega }$\emph{-distance} on the limit space $%
\mathcal{X}_{\omega }.$

The \emph{total distance} of the $\omega $-substitution graph system is
defined as%
\begin{equation*}
\mathcal{W}_{n}:=\sum_{\{u,v\}\subset \Sigma ^{n}}d_{n}(u,v),
\end{equation*}%
which can be viewed as a generalized Wiener index for a broader class of
metrics, encompassing the classical Wiener index and the Kirchhoff index as special cases.

Based on the GHP convergence established above, the asymptotic behavior of
the generalized Wiener index can be characterized by the continuous
integration of the limit metric. Since the normalized distance is a
continuous functional on the compact limit space, the leading-order
coefficient of $\mathcal{W}_{n}$ is precisely half of the total distance
integral over $\mathcal{X}_{\omega }.$

\begin{theorem}
\label{th-W}Let $\{(G_{n},d_{n})\}_{n\geq 1}$ be a metric-preserving $\omega 
$-substitution graph system. Then the normalized total distance satisfies
\begin{equation*}
\lim_{n\rightarrow \infty }\frac{\mathcal{W}_{n}}{2^{n}N^{2n}}=\frac{%
\mathcal{D}_{\omega }}{2}=\frac{(N-1)^{2}}{N(2N-1)}C_{\omega },
\end{equation*}%
where $C_{\omega }=2\sum_{k=2}^{\infty }2^{-k}d_{1}(\omega_{1},\omega_{k}).$
\end{theorem}

Furthermore, in Subsection \ref{subsection-closed} we derive closed-form
expressions for the generalized Wiener index; see Theorem \ref{th-W_n}.

The remainder of the paper is organized as follows. Section 3 collects the necessary preliminaries. Section 4 establishes the Gromov--Hausdorff--Prokhorov convergence of the associated compact spaces. Section 5 studies their fractal dimensions and prove Theorems \ref{th-dim}--\ref{th-dim-stab}. Section 6 is devoted to generalized Wiener indices. Section 7 concludes the paper.

\bigskip

\section{Preliminaries}

\subsection{Gromov--Hausdorff--Prokhorov metric}

\label{subsection-GHP}

Let $(X,d^{X})$ be a Polish metric space (i.e., complete and separable). For
two Borel subsets $A,B\in \mathcal{B}(X),$ the Hausdorff metric between $A$ and $%
B $ is 
\begin{equation*}
d_{\mathrm{H}}^{X}(A,B)=\inf \{\varepsilon >0:A\subset B^{\varepsilon }\text{
and }B\subset A^{\varepsilon }\},
\end{equation*}%
where $A^{\varepsilon }=\{x\in X:\inf_{y\in A}d^{X}(x,y)<\varepsilon \}.$
Let $\mathcal{M}_{f}(X)$ denote the set of all finite Borel measures on $X.$
If $\mu ,\nu \in \mathcal{M}_{f}(X),$ the Prokhorov metric between $\mu $
and $\nu $ is 
\begin{equation*}
d_{\mathrm{P}}^{X}(\mu ,\nu )=\inf \{\varepsilon >0:\mu (A)\leq \nu
(A^{\varepsilon })+\varepsilon \text{ and }\nu (A)\leq \mu (A^{\varepsilon
})+\varepsilon, \forall \text{ closed set }A\}.
\end{equation*}

Recall that a rooted weighted metric space $\mathcal{X}=(X,d,\emptyset ,\mu
) $ is a metric space $(X,d)$ with a distinguished element $\emptyset ,$
called the root, and a locally finite Borel measure $\mu .$ Two rooted
weighted metric spaces $\mathcal{X}=(X,d,\emptyset ,\mu )$ and $\mathcal{X}%
^{\prime }=(X^{\prime },d^{\prime },\emptyset ^{\prime },\mu ^{\prime })$
are said to be GHP-isometric if there exists a one-to-one isometric map $%
\Phi :X\rightarrow X^{\prime }$ such that $\Phi (\emptyset )=\emptyset
^{\prime }$ and $\Phi _{\ast }\mu =\mu ^{\prime }.$ In this case, $\Phi $ is
called a GHP-isometry. We now recall the definition of Gromov--Hausdorff--Prokhorov metric for compact rooted weighted metric
spaces. Let $\mathcal{X}=(X,d,\emptyset ,\mu )$ and $\mathcal{X}^{\prime
}=(X^{\prime },d^{\prime },\emptyset ^{\prime },\mu ^{\prime })$ be two such
spaces, and define%
\begin{equation}
d_{\mathrm{GHP}}(\mathcal{X},\mathcal{X}^{\prime }):=\inf_{\Phi ,\Phi
^{\prime },Z}(d^{Z}(\Phi (\emptyset ),\Phi ^{\prime }(\emptyset ^{\prime
}))+d_{\mathrm{H}}^{Z}(\Phi (X),\Phi ^{\prime }(X^{\prime }))+d_{\mathrm{P}%
}^{Z}(\Phi _{\ast }\mu ,\Phi _{\ast }^{\prime }\mu ^{\prime })),  \label{GHP}
\end{equation}%
where the infimum is taken over all isometric embeddings $\Phi
:X\hookrightarrow Z$ and $\Phi ^{\prime }:X^{\prime }\hookrightarrow Z$ into
a common Polish metric space $(Z,d^{Z}).$ Note that equation (\ref{GHP})
does not actually define a metric, since $d_{\mathrm{GHP}}(\mathcal{X},%
\mathcal{X}^{\prime })=0$ if $\mathcal{X}\ $and $\mathcal{X}^{\prime }\ $are
GHP-isometric. We therefore take $\mathbb{K}$ to be the set of GHP-isometry
classes of compact rooted weighted metric spaces, and identify each space
with its class.

It has been shown in \cite{Abraham} that $d_{\mathrm{GHP}}$ defines a metric
on the set $\mathbb{K}$. Furthermore, the metric space $(\mathbb{K},d_{%
\mathrm{GHP}})$ is a Polish space.

\subsection{Basic properties of $\protect\omega $-substitution graph systems}

\label{section-props}We denote by $\sigma $ the \emph{shift operator} on $%
\Sigma ^{\ast }\cup \Sigma ^{\mathbb{N}},$ defined by $\sigma
(x)=x_{2}x_{3}\cdots $ for $x=x_{1}x_{2}x_{3}\cdots .$ In particular, if $x$
is a finite word of length $n\geq 1,$ then $\sigma (x)$ is a word of length $%
n-1;$ for a single-letter word $x$, $\sigma (x)=\epsilon $.

The self-similarity of the sequence $\{G_{n}\}$ of $\omega $-substitution
graphs is naturally captured by the shift operator.

\begin{claim}
For any integers $n>1$ and $i\in \Sigma ,$ the restriction of the shift
operator $\sigma $ to the subset $i\Sigma ^{n-1}$ induces a graph
isomorphism from the induced subgraph $G_{n,i}:=G_{n}[i\Sigma ^{n-1}]$ to $%
G_{n-1}.$
\end{claim}

\begin{proof}
Let $u,v\in \Sigma ^{n-1}$. By the recursive construction of the $\omega $%
-substitution graphs, two vertices $iu$ and $iv$ in $G_{n}$ are adjacent if
and only if $u$ and $v$ are adjacent in $G_{n-1}$. Since the restriction $%
\sigma |_{i\Sigma ^{n-1}}:iu\mapsto u$ is a bijection between the vertex
sets $i\Sigma ^{n-1}$ and $\Sigma ^{n-1}$, and it preserves the adjacency
relation by the substitution rule, it follows that $\sigma |_{i\Sigma
^{n-1}} $ is a graph isomorphism.
\end{proof}

Let $\{(G_{n},d_{n})\}_{n\geq 1}$ be a metric-preserving $\omega $%
-substitution graph system. For each $n>1$ and $i\in \Sigma ,$ the set $%
\{i(\omega |_{n-1}),G_{n.i}\backslash \{i(\omega |_{n-1})\},G_{n}\backslash
G_{n.i}\}$ is a cutting of $G_{n},$ that is, the point $i\omega _{n-1}$ is a
cut-point of $G_{n},$ subsets $G_{n.i}\backslash \{i(\omega |_{n-1})\}$ and $%
G_{n}\backslash G_{n.i}$ form a separation of $G_{n}\backslash \{i(\omega
|_{n-1})\}.$ The following assertions follow from the definition and the self-similarity of the $\omega $-substitution graph system.

\begin{claim}
\label{clm-2}For any $u,v\in \Sigma ^{n}$ with $\left\vert u\wedge
v\right\vert =m,$ it holds
\begin{equation*}
d_{n}(u,v)=d_{n-t}(\sigma ^{t}u,\sigma ^{t}v)\text{ }\forall t\in \{1,\ldots
,m\}.
\end{equation*}
\end{claim}

\begin{claim}
\label{clm-3}For any $u,v\in \Sigma ^{n}$ with $\left\vert u\wedge
v\right\vert =0,$ it holds
\begin{equation*}
d_{n}(u,v)=d_{n-1}(\sigma (u),\omega
|_{n-1})+d_{1}(u_{1},v_{1})+d_{n-1}(\omega |_{n-1},\sigma (v)).
\end{equation*}
\end{claim}

Combining Claims \ref{clm-2} and \ref{clm-3}, one can get the following claim.

\begin{claim}
\label{clm-4}For any $u,v\in \Sigma ^{n}$ with $\left\vert u\wedge
v\right\vert =m<n,$ it holds
\begin{eqnarray*}
d_{n}(u,v) &=&d_{n-m-1}(\sigma ^{m+1}(u),\omega
|_{n-m-1})+d_{1}(u_{m+1},v_{m+1}) \\
&&+d_{n-m-1}(\omega |_{n-m-1},\sigma ^{m+1}(v)).
\end{eqnarray*}
\end{claim}

For $n>1$ and $1\leq k<n,$ we denote%
\begin{equation*}
\delta _{k,n}:=d_{k}(\omega |_{k},\sigma ^{n-k}(\omega |_{n})).
\end{equation*}%
The lemma below provides its recursive structure as well as its general
solution.

\begin{lemma}
\label{le-delta}For each $1<k\leq n-1,$ 
\begin{equation}
\delta _{k,n}=\sum_{l=1}^{k-1}\delta _{l,l+1}+\sum_{l=0}^{k-1}D_{n-l},
\label{delta_k,n}
\end{equation}%
where $D_{j}=d_{1}(\omega _{1},\omega _{j})$ for each $j\geq 2.$ In particular, 
\begin{equation}
\delta _{k,k+1}=2^{k+1}\sum_{p=2}^{k+1}2^{-p}D_{p}.  \label{delta_l,l+1}
\end{equation}%
Therefore, 
\begin{equation}
\delta
_{k,n}=2^{k+1}\sum_{p=2}^{k}2^{-p}D_{p}-\sum_{p=2}^{k}D_{p}+%
\sum_{p=n-k+1}^{n}D_{p}.  \label{delta_k,n-new}
\end{equation}
\end{lemma}

\begin{proof}
Since (\ref{delta_l,l+1})--(\ref{delta_k,n-new}) follow directly from (\ref%
{delta_k,n}), it suffices to prove (\ref{delta_k,n}). It follows from the
assumption (\ref{assump-omega}) and Claim \ref{clm-3} that 
\begin{eqnarray}
\delta _{k,n} &=&d_{k-1}(\sigma (\omega |_{k}),\omega |_{k-1})+d_{1}(\omega
_{1},\omega _{n-k+1})+d_{k-1}(\omega |_{k-1},\sigma ^{n-k+1}(\omega |_{n})) 
\notag \\
&=&\delta _{k-1,k}+d_{1}(\omega_{1},\omega _{n-k+1})+\delta _{k-1,n}.
\label{delta_k,n_recu}
\end{eqnarray}%
Note that $\delta _{1,n}=d_{1}(\omega _{1},\sigma ^{n-1}(\omega|
_{n}))=D_{n}, $ the recursive equation (\ref{delta_k,n_recu}) gives us (\ref{delta_k,n}).
\end{proof}

\begin{corollary}
\label{coro-conv}As $k\rightarrow \infty ,$ the ratio $\delta _{k,n}/2^{k}$
converges to $C_{\omega }$ uniformly in $n>k,$ where $C_{\omega
}=2\sum_{p=2}^{\infty }2^{-p}D_{p}$ is a constant depending only on $ \omega $.
\end{corollary}

For each $n\geq 1,$ we denote by $\mathrm{diam}(G_{n},d_{n})$ the diameter
of the graph $G_{n}$ with respect to the metric $d_{n}$, defined as 
\begin{equation*}
\mathrm{diam}(G_{n},d_{n})=\max_{u,v\in \Sigma ^{n}}d_{n}(u,v).
\end{equation*}%
For the sake of brevity, we write it as $\mathrm{diam}(G_{n})$, whenever the
metric is clear from the context.

\begin{proposition}
\label{prop-diam}Let $\{(G_{n},d_{n})\}$ be a metric-preserving $\omega $%
-substitution graph system. Then $\mathrm{diam}(G_{n})\asymp 2^{n}.$
\end{proposition}

\begin{proof}
It follows from (\ref{delta_l,l+1}) that for each $n\geq 2$, 
\begin{equation*}
\delta _{n-1,n}\geq 2^{n}2^{-2}D_{2},
\end{equation*}%
which implies 
\begin{equation*}
\mathrm{diam}(G_{n})\geq \delta _{n-1,n}\gtrsim 2^{n}.
\end{equation*}

On the other hand, from the self-similarity of $\{G_{n}\}$ and Claim \ref%
{clm-3}$,$ we have%
\begin{equation*}
\mathrm{diam}(G_{n})\leq 2\mathrm{diam}(G_{n-1})+D_{\max },
\end{equation*}%
where $D_{\max }=\max_{i,j\in \Sigma }d_{1}(i,j).$ A simple induction on $n$
then yields $\mathrm{diam}(G_{n})\lesssim 2^{n},$ which completes the proof.
\end{proof}

\bigskip

\section{Scaling limit of substitution graph systems}

In this section, let $\{(G_{n},d_{n})\}_{n\geq 1}$ be a metric-preserving $%
\omega $-substitution graph system.

For $x\in\Sigma^{\mathbb N}$, define the (possibly infinite) indices of
disagreement $\{\ell_p(x)\}_{p\ge1}$ between $x$ and $\omega$ recursively by
\begin{equation*}
\ell_0(x):=0,\quad
\ell_p(x):=\min\{k>\ell_{p-1}(x):x_k\neq\omega_k\},
\quad p\ge1,
\end{equation*}%
with the conventions $\min\varnothing:=\infty$ and $2^{-\infty}:=0$.
It follows from\ Claim \ref{clm-4} that%
\begin{equation}
d_{n}(x|_{n},\omega |_{n})=\sum_{\ell _{p}(x)<n}\delta _{n-\ell
_{p}(x),n}+\sum_{\ell _{p}(x)\leq n}d_{1}(x_{\ell _{p}(x)},\omega _{\ell
_{p}(x)}).  \label{d_n(x|n}
\end{equation}

For each $n\in \mathbb{N}$ and $u\in \Sigma ^{n},$ let $h_{n}(u):=\widetilde{%
d}_{n}(u,\omega |_{n})$.

\begin{lemma}
For any $x\in \Sigma ^{\mathbb{N}},$ the limit $h(x):=\lim\limits_{n\rightarrow\infty}h_{n}(x|_{n})$ exists. \newline
Moreover, one has
\begin{equation}
h(x)=C_{\omega }\sum_{p=1}^{\infty }2^{-\ell _{p}(x)}\text{ for }x\in \Sigma
^{\mathbb{N}}.  \label{h(x)}
\end{equation}
\end{lemma}

\begin{proof}
Given $x\in \Sigma ^{\mathbb{N}},$ from (\ref{d_n(x|n}), we only need to
consider 
\begin{equation*}
a_{n}(x):=2^{-n}\sum_{\ell _{p}(x)<n}\delta _{n-\ell _{p}(x),n}.
\end{equation*}%
Note that 
\begin{equation*}
a_{n}(x)=\sum_{p=1}^{\infty }a_{n,p}(x),\text{ where }a_{n,p}(x)=\mathbf{1}%
_{\{\ell _{p}(x)<n\}}\frac{\delta _{n-\ell _{p}(x),n}}{2^{n-\ell _{p}(x)}}%
2^{-\ell _{p}(x)}.
\end{equation*}%
Since $\delta _{k,n}/2^{k}\rightarrow C_{\omega }$ $(k<n),$ the family $%
\{\delta _{k,n}/2^{k}:1\leq k<n\}$ is uniformly bounded. Hence $\left\vert
a_{n,p}(x)\right\vert \leq A2^{-\ell _{p}(x)}$ for certain constant $A<\infty ,$
and $\sum_{p=1}^{\infty }2^{-\ell _{p}(x)}<\infty .$ Therefore the Dominated
Convergence Theorem yields%
\begin{equation*}
\lim_{n\rightarrow \infty }h_{n}(x|_n)=\sum_{p=1}^{\infty }\lim_{n\rightarrow
\infty }a_{n,p}(x)=C_{\omega }\sum_{p=1}^{\infty }2^{-\ell _{p}(x)}.
\end{equation*}
\end{proof}

Since $\sum_{p=1}^{\infty }2^{-\ell _{p}(x)}\leq 1$ for all $x\in \Sigma ^{\mathbb{N}}$, we have the following conclusion.

\begin{corollary}
\label{coro-h}The function $h(x)$ is bounded above by $C_{\omega }\ $and
vanishes if and only if $x=\omega .$
\end{corollary}

For any distinct finite words $u,v\in \Sigma _{n}$ with $\left\vert u\wedge
v\right\vert =m<n$, by Claims \ref{clm-4} 
\begin{equation}
\widetilde{d}_{n}(u,v)=2^{-m-1}\left( h_{n-m-1}(\sigma
^{m+1}u)+h_{n-m-1}(\sigma ^{m+1}v)\right) +2^{-n}d_{1}(u_{m+1},v_{m+1}).
\label{d_n=2^-}
\end{equation}%
This leads to the following proposition.

\begin{proposition}
\label{prop-d-limit}For any $x,y\in \Sigma ^{\mathbb{N}},$ the limit $%
\widetilde{d}(x,y):=\lim\limits_{n\rightarrow\infty}\widetilde{d}%
_{n}(x|_{n},y|_{n})$ exists.\newline
Moreover, for any distinct $x,y\in \Sigma ^{\mathbb{N}},$
\begin{equation}
\widetilde{d}(x,y)=2^{-\left\vert x\wedge y\right\vert -1}(h(\sigma
^{\left\vert x\wedge y\right\vert +1}x)+h(\sigma ^{\left\vert x\wedge
y\right\vert +1}y)).  \label{d=2^-(h+h)}
\end{equation}
\end{proposition}

\begin{proof}
It suffices to prove the case $x\neq y.$ Let $m=\left\vert x\wedge
y\right\vert .$ By Claims \ref{clm-4}, we have 
\begin{eqnarray*}
\widetilde{d}_{n}(x|_{n},y|_{n}) &=&2^{-m-1}\left( h_{n-m-1}(\sigma
^{m+1}(x|_{n}))+h_{n-m-1}(\sigma ^{m+1}(y|_{n}))\right) \\
&&+2^{-n}d_{1}(x_{m+1},y_{m+1}) \\
&=&2^{-m-1}\left( h_{n-m-1}((\sigma ^{m+1}x)|_{n-m-1})+h_{n-m-1}((\sigma
^{m+1}y)|_{n-m-1})\right) \\
&&+2^{-n}d_{1}(x_{m+1},y_{m+1}).
\end{eqnarray*}%
As $n\rightarrow \infty ,$ it follows that 
\begin{equation*}
\lim_{n\rightarrow \infty }\widetilde{d}_{n}(x|_{n},y|_{n})=2^{-m-1}(h(%
\sigma ^{m+1}x)+h(\sigma ^{m+1}y)).
\end{equation*}
\end{proof}

From Corollary \ref{coro-h}, we obtain the following result.

\begin{corollary}
\label{coro-d}The function $\widetilde{d}(x,y)$ is bounded above by $%
2^{-\left\vert x\wedge y\right\vert }C_{\omega }$ and vanishes if and only
if 
\begin{equation*}
\sigma ^{\left\vert x\wedge y\right\vert +1}x=\omega =\sigma ^{\left\vert x\wedge
y\right\vert +1}y \text{ or } x=y.
\end{equation*}
\end{corollary}

By the previous results, $\widetilde{d}$ is a pseudometric on $\Sigma ^{%
\mathbb{N}}.$ Define an equivalence relation on $\Sigma ^{\mathbb{N}}$ by 
\begin{equation*}
x\sim y\Longleftrightarrow \widetilde{d}(x,y)=0.
\end{equation*}%
Let $X_{\omega }$ be the corresponding quotient space, and $\pi :\Sigma ^{%
\mathbb{N}}\rightarrow X_{\omega }$ be the canonical projection mapping each
infinite word $x$ to its equivalence class $\left\langle x\right\rangle .$
Equip $X_{\omega }$ with the quotient metric 
\begin{equation*}
d_{\omega }(\left\langle x\right\rangle ,\left\langle y\right\rangle )=%
\widetilde{d}(x,y).
\end{equation*}%
Since $\widetilde{d}$ is a pseudometric, the above definition is independent
of the choice of representatives and defines a metric on $X_{\omega }.$
Since $\Sigma ^{\mathbb{N}}$ is compact, the quotient space $X_{\omega }$ is
compact. Hence $(X_{\omega },d_{\omega })$ is a compact metric space.

Recall $\mu $ is the uniform Bernoulli measure on $\Sigma ^{\mathbb{N}}.$ Define a Borel probability measure $\mu _{\omega }$ on $X_{\omega }$ as the
pushforward of $\mu $ under $\pi ,$ that is, 
\begin{equation*}
\mu _{\omega }:=\pi _{\sharp }\mu ,\text{ where }\mu _{\omega }(A)=\mu (\pi
^{-1}(A))
\end{equation*}%
for every Borel set $A\subset X_{\omega }.$ Since $\pi $ is measurable and $%
\mu $ is a probability measure on $\Sigma ^{\mathbb{N}},$ the measure $\mu
_{\omega }$ is well-defined and satisfies $\mu _{\omega }(X_{\omega })=1.$

Let $\mathcal{X}_{\omega }=(X_{\omega },d_{\omega },\left\langle \overline{1}%
\right\rangle ,\mu _{\omega }).$ Both $\mathcal{X}_{n}^{\omega }$ and $%
\mathcal{X}_{\omega }$ are compact rooted weighted metric spaces.
Identifying compact rooted weighted metric spaces with their GHP-isometry
classes, we regard them as elements of $\mathbb{K}.$ We now provide a precise statement of the main convergence theorem.

\begin{theorem}
\label{th-GHP-1}$\mathcal{X}_{n}^{\omega }$ converges to $\mathcal{X}%
_{\omega }$ in the Gromov--Hausdorff--Prokhorov topology.
\end{theorem}

For a cylinder set $[u],$ let $[u]_{\omega }:=\pi ([u])\subset X_{\omega },$
which induces a natural family of projected cylinder sets in $X_{\omega }.$

While the cylinder sets $\{[u]:u\in \Sigma ^{n}\}$ form a partition of $%
\Sigma ^{\mathbb{N}},$ their images under $\pi $ may overlap in the quotient
space $X_{\omega }.$ The following lemma shows that these projected sets
remain disjoint in a measure-theoretic sense. Consequently, for any finite word $u\in \Sigma ^{\ast },$ one has
\begin{equation*}
\mu _{\omega }([u]_{\omega })=\mu ([u])=N^{-\left\vert u\right\vert }.
\end{equation*}%

\begin{lemma}
\label{le-disjoint-mu}For any distinct $u,v\in \Sigma ^{n},$ it holds that%
\begin{equation*}
\mu _{\omega }([u]_{\omega }\cap \lbrack v]_{\omega })=0.
\end{equation*}
\end{lemma}

\begin{proof}
Suppose that $\left\vert u\wedge v\right\vert =m<n.$ Let $\xi \in \lbrack
u]_{\omega }\cap \lbrack v]_{\omega },$ then there exist $z_{1}\in \lbrack
u] $ and $z_{2}\in \lbrack v]$ such that $\xi =\left\langle z_{1}\right\rangle=\left\langle z_{2}\right\rangle.$ Hence $%
z_{1}\sim z_{2}.$ By Corollary \ref{coro-d}, we have $z_{1}=(u|_{m})u_{m+1}%
\omega $ and $z_{2}=(u|_{m})v_{m+1}\omega ,$ where $u_{m+1}\neq v_{m+1}.$

Let $x\in \pi ^{-1}(\xi ).$ Since $\omega $ is not periodic and $x\sim z_{1}$, Corollary \ref{coro-d} again implies $x=(u|_{m})q\omega $
for certain $q\in \Sigma $. Therefore%
\begin{equation*}
\pi ^{-1}([u]_{\omega }\cap \lbrack v]_{\omega })\subset \{z=u|_{m}q\omega
:q\in \Sigma \}.
\end{equation*}%
The latter set is finite, hence $\mu (\pi ^{-1}([u]_{\omega }\cap \lbrack
v]_{\omega }))=0.$ By the definition of the push-forward measure, we have 
\begin{equation*}
\mu _{\omega }([u]_{\omega }\cap \lbrack v]_{\omega })=\mu (\pi
^{-1}([u]_{\omega }\cap \lbrack v]_{\omega }))=0.
\end{equation*}
\end{proof}

For each $n\in \mathbb{N},$ define $\iota _{n}:\Sigma ^{n}\rightarrow
X_{\omega }$ as 
\begin{equation*}
\iota _{n}(v)=\left\langle v\overline{1}\right\rangle ,\text{ }v\in \Sigma
^{n}.
\end{equation*}

\begin{lemma}
\label{le-dis_iota}The sequence 
\begin{equation*}
\mathrm{dis}(\iota _{n}):=\sup_{u,v\in \Sigma ^{n}}\left\vert \widetilde{d}%
_{n}(u,v)-d_{\omega }(\iota _{n}(u),\iota _{n}(v))\right\vert
\end{equation*}%
converges to zero as $n\rightarrow \infty .$
\end{lemma}

\begin{proof}
For each $n\in \mathbb{N}$, let $u,v\in \Sigma ^{n}$ with $\left\vert
u\wedge v\right\vert =m<n.$ Combining (\ref{d_n=2^-})--(\ref{d=2^-(h+h)}),
\begin{eqnarray}
&&\left\vert d_{\omega }(\iota _{n}(u),\iota _{n}(v))-\widetilde{d}%
_{n}(u,v)\right\vert =\left\vert \widetilde{d}(u\overline{1},v\overline{1})-%
\widetilde{d}_{n}(u,v)\right\vert  \notag \\
&\leq &2^{-m-1}\left\vert h(\sigma ^{m+1}\left( u\overline{1}\right)
)-h_{n-m-1}(\sigma ^{m+1}u)\right\vert  \label{d_omega-d_n} \\
&&+2^{-m-1}\left\vert h(\sigma ^{m+1}\left( v\overline{1}\right)
)-h_{n-m-1}(\sigma ^{m+1}v)\right\vert +2^{-n}d_{1}(u_{m+1},v_{m+1}).  \notag
\end{eqnarray}

Write $n^{\prime }=n-m-1$ and $x^{\prime }=\sigma ^{m+1}\left( u\overline{1}%
\right) .$ It follows from (\ref{d_n(x|n})--(\ref{h(x)}) and (\ref%
{delta_k,n-new}) that%
\begin{equation}
\begin{aligned} &\left| h(\sigma^{m+1} (u\overline{1})) -
h_{n-m-1}(\sigma^{m+1}u) \right| = \left| h(x') - h_{n-m-1}(x'|_{n-m-1})
\right| \\ & \le \left| C_{\omega} \sum_{p=1}^{\infty} 2^{-\ell_p(x')}
-  \sum_{\ell_p(x') < n'} 2^{\ell_p(x')+1} \sum_{j=2}^{n'-\ell_p} 2^{-j} D_j
\right| \\ &\quad + 2^{-n'} \Biggl( \sum_{\ell_p(x') < n'} \biggl(
\sum_{j=2}^{n'-\ell_p(x')} D_j + \sum_{j=\ell_p(x')+1}^{n'} D_j \biggr) +
\sum_{\ell_p(x') \le n'} d_1(x_{\ell_p(x')}, \omega_{\ell_p(x')}) \Biggr).
\end{aligned}  \label{1234567}
\end{equation}%
The first and second terms on the right-hand side of (\ref{1234567}), denoted by $\mathbf{I}_{1}$ and $\mathbf{I}_{2}$, respectively, admit the estimates%
\begin{eqnarray*}
\mathbf{I}_{1} &=&\sum_{\ell _{p}(x^{\prime })<n^{\prime }}2^{-\ell
_{p}(x^{\prime })}\left( C_{\omega }-2\sum_{j=2}^{n^{\prime }-\ell
_{p}}2^{-j}D_{j}\right) +C_{\omega }\sum_{\ell _{p}(x^{\prime })\geq
n^{\prime }}2^{-\ell _{p}(x^{\prime })} \\
&\leq &\sum_{\ell _{p}(x^{\prime })<n^{\prime }}2^{-\ell _{p}(x^{\prime
})}\left( 2\sum_{j=n^{\prime }-\ell _{p}(x^{\prime })+1}^{\infty
}2^{-j}D_{j}\right) +C_{\omega }2^{-n^{\prime }+1} \\
&\leq &\sum_{\ell _{p}(x^{\prime })<n^{\prime }}2^{-\ell _{p}(x^{\prime
})}\left( 2D_{\max }2^{-n^{\prime }+\ell _{p}(x^{\prime })}\right)
+C_{\omega }2^{-n^{\prime }+1} \\
&\leq &2^{-n^{\prime }}(2nD_{\max }+2C_{\omega }),
\end{eqnarray*}%
and 
\begin{equation*}
\mathbf{I}_{2}\leq 2^{-n^{\prime }}(2(n^{\prime })^{2}+n^{\prime })D_{\max
}\leq 2^{-n^{\prime }}(2n^{2}+n)D_{\max }.
\end{equation*}%
Combining these results, it follows that%
\begin{equation}
\left\vert h(\sigma ^{m+1}( u\overline{1}) )-h_{n-m-1}(\sigma
^{m+1}u)\right\vert \leq 2^{-n^{\prime }}P(n),  \label{123}
\end{equation}%
where $P(n)=2C_{\omega }+(2n^{2}+3n)D_{\max }.$ Similarly, the same estimate
holds for the term involving $v,$ yielding%
\begin{equation}
\left\vert h(\sigma ^{m+1}\left( v\overline{1}\right) )-h_{n-m-1}(\sigma
^{m+1}v)\right\vert \leq 2^{-n^{\prime }}P(n).  \label{1234}
\end{equation}%
Substituting (\ref{123})--(\ref{1234}) into (\ref{d_omega-d_n}) yields 
\begin{equation*}
\left\vert d_{\omega }(\iota _{n}(u),\iota _{n}(v))-\widetilde{d}%
_{n}(u,v)\right\vert \leq 2^{-n}(2P(n)+D_{\max }),
\end{equation*}%
which implies 
\begin{equation*}
\mathrm{dis}(\iota _{n})\leq 2^{-n}(2P(n)+D_{\max }) \to 0 \quad (n\to\infty).
\end{equation*}
\end{proof}

To study the GHP distance, we now consider the disjoint union $Z_{n}=\Sigma
^{n}\sqcup X_{\omega }$. We shall abuse notation by identifying $\Sigma ^{n}$
and $X_{\omega }$ (along with their respective metrics and measures) with
their images under the canonical embeddings into $Z_{n}.$

Define a function $d_{Z_{n}}$ on $Z_{n}\times Z_{n}$ by 
\begin{equation*}
d_{Z_{n}}(x,y):=\left\{ 
\begin{array}{ll}
\widetilde{d}_{n}(x,y), & \text{if }x,y\in \Sigma ^{n}, \\ 
d_{\omega }(x,y), & \text{if }x,y\in X_{\omega }, \\ 
\inf\limits_{u\in \Sigma ^{n}}\{\widetilde{d}_{n}(x,u)+\varepsilon
_{n}+d_{\omega }(\iota _{n}(u),y)\}, & \text{otherwise,}%
\end{array}%
\right.
\end{equation*}%
where the last case is considered up to symmetry and $\varepsilon _{n}:=\frac{%
1}{2}\mathrm{dis}(\iota _{n})+2^{-n}.$

\begin{lemma}
\label{le-dZn}The function $d_{Z_{n}}$ is a metric on $Z_{n}.$
\end{lemma}

\begin{proof}
It is clear from the definition that $d_{Z_{n}}$ is non-negative, symmetric,
and $d_{Z_{n}}(x,y)=0$ if and only if $x=y.$ It remains to verify the
triangle inequality $d_{Z_{n}}(x,y)\leq d_{Z_{n}}(x,z)+d_{Z_{n}}(z,y)$ for any $%
x,y,z\in Z_{n}.$

If all three points lie in $\Sigma ^{n}$ or $X_{\omega },$ the
triangle inequality holds as $\widetilde{d}_{n}$ and $d_{\omega }$ are
metrics. By symmetry, we only need to consider the following several cases.

\textbf{Case 1.} Let $x,z\in \Sigma ^{n}$ and $y\in X_{\omega },$ then the
triangle inequality for $\widetilde{d}_{n}$ implies 
\begin{eqnarray*}
d_{Z_{n}}(x,y) &\leq &\inf_{u\in \Sigma _{n}}\{(\widetilde{d}_{n}(x,z)+%
\widetilde{d}_{n}(z,u))+\varepsilon _{n}+\widetilde{d}(\iota _{n}(u),y)\} \\
&=&\widetilde{d}_{n}(x,z)+d_{Z_{n}}(z,y)=d_{Z_{n}}(x,z)+d_{Z_{n}}(z,y).
\end{eqnarray*}%
The case where $x,z\in X_{\omega }$ and $y\in \Sigma ^{n}$ follows by an
analogous argument using the triangle inequality for $d_{\omega }.$

\textbf{Case 2}. Let $x,y\in \Sigma ^{n}$ and $z\in X_{\omega }.$ For any $%
u^{\prime },v^{\prime }\in \Sigma _{n},$ combining the triangle inequality
for $d_{\omega }$ and $\widetilde{d}_{n}$ with the definition of $%
\varepsilon _{n},$ it follows that%
\begin{eqnarray*}
&&\widetilde{d}_{n}(x,u^{\prime })+\varepsilon _{n}+\widetilde{d}(\iota
_{n}(u^{\prime }),z)+\widetilde{d}_{n}(y,v^{\prime })+\varepsilon _{n}+%
\widetilde{d}(\iota _{n}(v^{\prime }),z) \\
&\geq &\widetilde{d}_{n}(x,u^{\prime })+\widetilde{d}_{n}(y,v^{\prime
})+2\varepsilon _{n}+\widetilde{d}(\iota _{n}(u^{\prime }),\iota
_{n}(v^{\prime })) \\
&\geq &\widetilde{d}_{n}(x,u^{\prime })+\widetilde{d}_{n}(y,v^{\prime })+%
\widetilde{d}_{n}(u^{\prime },v^{\prime }) \\
&\geq &\widetilde{d}_{n}(x,y)=d_{Z_{n}}(x,y).
\end{eqnarray*}%
Taking the infimum over $u^{\prime }$ and $v^{\prime },$ we obtain $%
d_{Z_{n}}(x,z)+d_{Z_{n}}(z,y)\geq d_{Z_{n}}(x,y).$

\textbf{Case 3.} Let $x,y\in X_{\omega }$ and $z\in \Sigma ^{n}.$ Similar to
Case 2, for any $u^{\prime },v^{\prime }\in \Sigma ^{n},$ 
\begin{eqnarray*}
&&\widetilde{d}_{n}(z,u^{\prime })+\varepsilon _{n}+\widetilde{d}(\iota
_{n}(u^{\prime }),x)+\widetilde{d}_{n}(z,v^{\prime })+\varepsilon _{n}+%
\widetilde{d}(\iota _{n}(v^{\prime }),y) \\
&\geq &\widetilde{d}_{n}(u^{\prime },v^{\prime })+2\varepsilon _{n}+%
\widetilde{d}(\iota _{n}(u^{\prime }),x)+\widetilde{d}(\iota _{n}(v^{\prime
}),y) \\
&\geq &\widetilde{d}(\iota _{n}(u^{\prime }),\iota _{n}(v^{\prime }))+%
\widetilde{d}(\iota _{n}(u^{\prime }),x)+\widetilde{d}(\iota _{n}(v^{\prime
}),y) \\
&\geq &\widetilde{d}(x,y)=d_{Z_{n}}(x,y).
\end{eqnarray*}%
The conclusion follows by taking the infimum over $u^{\prime }$ and $%
v^{\prime }.$
\end{proof}

The Hausdorff distance admits the following estimate.
\begin{proposition}
\label{prop-d_H}For each $n\in \mathbb{N},$ the Hausdorff distance between $%
\Sigma ^{n}$ and $X_{\omega }$ in $(Z_{n},d_{Z_{n}})$ satisfies%
\begin{equation*}
d_{\mathrm{H}}^{d_{Z_{n}}}(\Sigma ^{n},X_{\omega })< \varepsilon
_{n}^{\prime },
\end{equation*}%
where $\varepsilon _{n}^{\prime }=\varepsilon _{n}+2^{-n}(C_{\omega }+1).$
\end{proposition}

\begin{proof}
For any $ \xi =\left\langle x\right\rangle \in X_{\omega },$ applying
Corollary \ref{coro-d} we have 
\begin{eqnarray*}
d_{Z_{n}}(x|_{n},\xi ) &\leq &\widetilde{d}_{n}(x|_{n},x|_{n})+\varepsilon
_{n}+d_{\omega }(\iota _{n}(x|_{n}),\xi ) \\
&=&\varepsilon _{n}+\widetilde{d}((x|_{n})\overline{1},x) \\
&\leq &\varepsilon _{n}+2^{-n}C_{\omega }<\varepsilon _{n}^{\prime },
\end{eqnarray*}%
hence $d_{Z_{n}}(\xi ,\Sigma ^{n})<\varepsilon _{n}^{\prime }.$

On the other hand, for each $v\in \Sigma ^{n},$ it follows that 
\begin{equation*}
d_{Z_{n}}(v,\iota _{n}(v))\leq \widetilde{d}_{n}(v,v)+\varepsilon
_{n}+d_{\omega }(\iota _{n}(v),\iota _{n}(v))=\varepsilon _{n},
\end{equation*}%
which implies $d_{Z_{n}}(v,X_{\omega })\leq \varepsilon _{n}.$

Combining these two estimates, we conclude that $d_{\mathrm{H}%
}^{d_{Z_{n}}}(\Sigma ^{n},X_{\omega })< \varepsilon _{n}^{\prime }.$
\end{proof}

\begin{lemma}
\label{le-v}For any $v\in \Sigma ^{n},$ it holds that 
\begin{equation*}
\lbrack v]_{\omega }\subset B_{d_{Z_{n}}}(v,\varepsilon _{n}^{\prime }).
\end{equation*}
\end{lemma}

\begin{proof}
Note that $x\in \lbrack v]_{\omega }$ if and only if there exists a $z\in
\lbrack v]$ such that $[z]=x.$

If $x\in \lbrack v]_{\omega },$ then 
\begin{equation*}
d_{\omega }(\iota _{n}(v),x)\leq \widetilde{d}(v\overline{1},z)=2^{-n-1}(h(%
\overline{1})+h(\sigma ^{n+1}z))\leq 2^{-n}C_{\omega },
\end{equation*}%
which implies 
\begin{equation*}
d_{Z_{n}}(v,x)\leq \widetilde{d}_{n}(v,v)+\varepsilon _{n}+d_{\omega }(\iota
_{n}(v),x)< \varepsilon _{n}^{\prime }.
\end{equation*}
\end{proof}

A corresponding estimate holds for the Prokhorov distance.
\begin{proposition}
\label{prop-d_P}For each $n\in \mathbb{N},$ the Prokhorov distance between $%
\mu _{n}$ and $\mu _{\omega }$ in $(Z_{n},d_{Z_{n}})$ satisfies%
\begin{equation*}
d_{\mathrm{P}}^{d_{Z_{n}}}(\mu _{n},\mu _{\omega })\leq \varepsilon
_{n}^{\prime }.
\end{equation*}
\end{proposition}

\begin{proof}
Suppose $A\subset Z_{n}$ is a Borel subset. If $A\subset X_{\omega },$ then
by Lemmas \ref{le-disjoint-mu} and \ref{le-v}%
\begin{eqnarray*}
\mu _{\omega }(A) &=&\mu _{\omega }(\bigcup_{v\in \Sigma _{n}}(A\cap \lbrack
v]_{\omega }))=\sum_{v\in \Sigma _{n}}\mu _{\omega }(A\cap \lbrack
v]_{\omega }) \\
&\leq &\sum_{v\in \Sigma _{n},A\cap \lbrack v]_{\omega }\neq \varnothing
}\mu _{\omega }([v]_{\omega })=\sum_{v\in \Sigma _{n},A\cap \lbrack
v]_{\omega }\neq \varnothing }\mu _{n}(\{v\}) \\
&\leq &\sum_{v\in \Sigma _{n}\cap A^{\varepsilon _{n}^{\prime }}}\mu
_{n}(\{v\})=\mu _{n}(A^{\varepsilon _{n}^{\prime }}).
\end{eqnarray*}%
Similarly, if $A\subset \Sigma _{n},$ then 
\begin{eqnarray*}
\mu _{n}(A) &=&\sum_{v\in A}\mu _{n}(\{v\})=\sum_{v\in A}\mu _{\omega
}([v]_{\omega }) =               \mu _{\omega }(\bigcup_{v\in A}[v]_{\omega })\leq \mu _{\omega
}(A^{\varepsilon _{n}^{\prime }}).
\end{eqnarray*}%
Combining the above two steps with the definition of $d_{P}^{d_{Z_{n}}},$ one has $d_{P}^{d_{Z_{n}}}(\mu _{n},\mu )\leq \varepsilon _{n}^{\prime }.$
\end{proof}

\begin{proof}[Proof of Theorem \protect\ref{th-GHP-1}]
According to Lemma \ref{le-dZn} and the definition of the GHP metric, 
\begin{equation*}
d_{\mathrm{GHP}}(\mathcal{X}_{n},\mathcal{X})\leq d_{Z_{n}}(\overline{1}%
|_{n},\overline{1})+d_{\mathrm{H}}^{d_{Z_{n}}}(\Sigma ^{n},X_{\omega })+d_{%
\mathrm{P}}^{d_{Z_{n}}}(\mu _{n},\mu _{\omega }).
\end{equation*}%
By the definition of $d_{Z_{n}},$ we obtain 
\begin{equation*}
d_{Z_{n}}(\overline{1}|_{n},\overline{1})\leq \widetilde{d}_{n}(\overline{1}%
|_{n},\overline{1}|_{n})+\varepsilon _{n}+\widetilde{d}(\overline{1},%
\overline{1})=\varepsilon _{n}.
\end{equation*}%
It follows from Propositions \ref{prop-d_H}--\ref{prop-d_P} and Lemma \ref%
{le-dis_iota} that 
\begin{equation*}
d_{\mathrm{GHP}}(\mathcal{X}_{n},\mathcal{X})\leq \varepsilon
_{n}+2\varepsilon _{n}^{\prime }\rightarrow 0\text{ as }n\rightarrow \infty .
\end{equation*}
\end{proof}

\bigskip

\section{Fractal Dimensions and Scaling Consistency}

This section studies the dimensional properties of both the limiting metric space $%
X_{\omega }$ and the $\omega $-substitution graph system. The Hausdorff and box-counting dimensions of $%
X_{\omega }$ as well as the discrete Hausdorff dimension of the graph sequence, will all be shown to equal $s:=\log N/\log 2.$ 

\begin{proof}[Proof of Theorem \protect\ref{th-dim}]
According to the theory of fractal dimensions, it suffices to prove that $%
\overline{\dim }_{\mathrm{B}}(X_{\omega })\leq s$ and $\dim _{\mathrm{H}%
}(X_{\omega })\geq s.$

For the upper bound, we consider the symbolic space $(\Sigma ^{\mathbb{N}%
},\rho )$ equipped with the standard metric $\rho (x,y)=2^{-\left\vert
x\wedge y\right\vert }.$ It is a well-known result in fractal geometry that
the Hausdorff and box-counting dimensions of $(\Sigma ^{\mathbb{N}},\rho )$
are both equal to $s.$ By Corollary \ref{coro-d}, the canonical quotient map 
$\pi :(\Sigma ^{\mathbb{N}},\rho )\rightarrow (X_{\omega },d_{\omega })$ is
Lipschitz continuous. Consequently, we obtain $\overline{\dim }_{\mathrm{B}%
}(X_{\omega })\leq \dim _{\mathrm{B}}(\Sigma ^{\mathbb{N}})=s.$

We next prove the lower bound. Fix $M\geq 3$, and set 
\begin{equation*}
A_{M}:=\{k\in \mathbb{N}:k\equiv 1\ \text{or}\ 2\pmod M\}.
\end{equation*}%
Define 
\begin{equation*}
F_{M}:=\left\{ x\in \Sigma ^{\mathbb{N}}:x_{k}=\omega _{1}\ \text{for every }%
k\in A_{M}\right\} .
\end{equation*}%
We first show that 
\begin{equation*}
\ell _{1}(\sigma ^{k}x)\leq M\text{\ for all }x\in F_{M},\ k\geq 0.
\end{equation*}%
Indeed, every interval of $M$ consecutive positive integers contains exactly
two elements of $A_{M}$. Thus, among the integers $k+1,k+2,\ldots ,k+M,$
there are two indices belonging to $A_{M}$. At least one of them can be
written as $k+j$ for some $2\leq j\leq M$. For this $j$, 
\begin{equation*}
(\sigma ^{k}x)_{j}=x_{k+j}=\omega _{1}.
\end{equation*}%
By assumption~(\ref{assump-omega}), one has $(\sigma ^{k}x)_{j}\neq \omega _{j}.$
Hence $\ell _{1}(\sigma ^{k}x)\leq M$.

By~(\ref{h(x)}), one has 
\begin{equation*}
C_{\omega }2^{-M}\leq h(\sigma ^{k}x)\leq C_{\omega }\text{\ for all }x\in
F_{M},\ k\geq 0.
\end{equation*}%
Let $x,y\in F_{M}$ be distinct and set $m=|x\wedge y|$. By~(\ref{d=2^-(h+h)}%
), it follows that 
\begin{equation*}
C_{\omega }2^{-M}2^{-m}\leq d_{\omega }(\pi (x),\pi (y))\leq C_{\omega
}2^{-m}.
\end{equation*}%
Equivalently, 
\begin{equation*}
C_{\omega }2^{-M}\rho (x,y)\leq d_{\omega }(\pi (x),\pi (y))\leq C_{\omega
}\rho (x,y).
\end{equation*}%
Note that $C_{\omega }>0$ by assumption (\ref{assump-omega}).
Thus 
\begin{equation*}
\pi |_{F_{M}}:F_{M}\longrightarrow \pi (F_{M})
\end{equation*}%
is a bi-Lipschitz bijection.

It remains to estimate the Hausdorff dimension of $F_{M}$. For $n\geq 1$,
let 
\begin{equation*}
q_{M}(n):=\#(\{1,\ldots ,n\}\setminus A_{M}).
\end{equation*}%
For each $k\in \mathbb{N}$, define a probability measure $\lambda _{k}$ on $%
\Sigma $ by 
\begin{equation*}
\lambda _{k}=\left\{ 
\begin{array}{ll}
\delta _{\omega _{1}}, & k\in A_{M}, \\ 
N^{-1}\sum_{a\in \Sigma }\delta _{a}, & k\notin A_{M}.%
\end{array}%
\right.
\end{equation*}%
Let $\nu _{M}:=\bigotimes_{k=1}^{\infty }\lambda _{k}$ be the corresponding
product probability measure on $\Sigma ^{\mathbb{N}}$. Clearly, $\nu
_{M}(F_{M})=1,$ and hence $\nu _{M}$ may be regarded as a probability
measure on $F_{M}$. Moreover, if $[u]\cap F_{M}\neq \varnothing $ and $|u|=n$%
, then 
\begin{equation*}
\nu _{M}([u]\cap F_{M})=N^{-q_{M}(n)}.
\end{equation*}

Since $\#(A_{M}\cap \{1,\ldots ,n\})\leq \frac{2n}{M}+2,$ we have $%
q_{M}(n)\geq \frac{M-2}{M}n-2.$ Set 
\begin{equation*}
s_{M}:=\frac{M-2}{M}\frac{\log N}{\log 2}.
\end{equation*}%
For $x\in F_{M}$ and $0<r<1,$ choose $n\geq 1$ such that $2^{-n}\leq
r<2^{-(n-1)}.$ If $y\in B_{\rho }(x,r)$ with $y\neq x $, then $2^{-|x\wedge y|}=\rho (x,y)<2^{-(n-1)}$,
and hence $|x\wedge y|\geq n.$ Therefore $y|_{n}=x|_{n}$, and thus 
\begin{equation*}
B_{\rho }(x,r)\cap F_{M}\subset \lbrack x|_{n}]\cap F_{M}.
\end{equation*}%
Since $\nu _{M}(F_{M})=1$, it follows that 
\begin{eqnarray*}
\nu _{M}(B_{\rho }(x,r)) &=&\nu _{M}(B_{\rho }(x,r)\cap F_{M})\leq \nu
_{M}([x|_{n}]\cap F_{M}) \\
&=&N^{-q_{M}(n)}\leq N^{2}N^{-\frac{M-2}{M}n}=N^{2}2^{-ns_{M}}\leq
N^{2}r^{s_{M}},
\end{eqnarray*}%
where the last inequality follows from $2^{-n}\leq r$. By the Mass
Distribution Principle, 
\begin{equation*}
\dim _{H}(F_{M})\geq s_{M}.
\end{equation*}%
Since $\pi |_{F_{M}}$ is bi-Lipschitz, 
\begin{equation*}
\dim _{H}(X_{\omega })\geq \dim _{H}(\pi (F_{M}))=\dim _{H}(F_{M})\geq s_{M}.
\end{equation*}%
Letting $M\rightarrow \infty $, we obtain 
\begin{equation*}
\dim _{H}(X_{\omega })\geq s.
\end{equation*}
Combining this with the upper bound, we conclude that
\begin{equation*}
  \dim_H(X_\omega)
=
\dim_B(X_\omega)
=s.
\end{equation*}
\end{proof}

Following \cite{Zeng3}, we recall the notion of the discrete Hausdorff dimension for graph sequences.

\begin{definition}
Let $\{G_{n}\}$ be a sequence of finite graphs with $\sharp
V(G_{n})\rightarrow \infty $ and $d_{n}$ be the geodesic distance on $G_{n}.$
The discrete Hausdorff dimension of $\{G_{n}\}$ is defined by 
\begin{equation*}
\dim _{\mathrm{DH}}(\{G_{n}\}):=\sup \left\{ \alpha \geq 0
:\varliminf\limits_{n\rightarrow \infty }\frac{1}{(\sharp V(G_{n}))^{2}}\sum 
_{\substack{ u,v\in V(G_{n})  \\ u\neq v}}\left( \frac{d_{n}(u,v)}{\mathrm{%
diam}(G_{n})}\right) ^{-\alpha }<\infty \right\} .
\end{equation*}
\end{definition}

We now turn to the proof of Theorem 3.

\begin{proof}[Proof of Theorem \protect\ref{th-dim-stab}]
By Theorem \ref{th-dim}, it suffices to show that $\dim _{\mathrm{DH}}(%
\mathcal{G})=s.$

For abbreviation, we continue to write $\widetilde{d}_{n}:=d_{n}/2^{n}.$ Let 
\begin{equation*}
I_{n}^{\alpha } := \sum_{\substack{ u,v\in \Sigma ^{n}  \\ u\neq v}} \left( 
\widetilde{d}_{n}(u,v)\right) ^{-\alpha } \quad\text{and}\quad J_{n}^{\alpha
} := \sum_{\substack{ u,v\in \Sigma ^{n}  \\ u\neq v}} \left( \frac{%
d_{n}(u,v)}{\mathrm{diam}(G_{n})} \right) ^{-\alpha }.
\end{equation*}
It follows from Proposition \ref{prop-diam} that 
\begin{equation}
I_{n}^{\alpha }\asymp J_{n}^{\alpha }.  \label{I=J}
\end{equation}

If $\alpha>s,$ by considering only the internal pairs within each sub-block
we have 
\begin{equation*}
I_{n}^{\alpha } \geq \sum_{j\in\Sigma} \sum_{\substack{ u,v\in j\Sigma^{n-1} 
\\ u\neq v}} \left(\widetilde d_n(u,v)\right)^{-\alpha} = 2^{\alpha}N
I_{n-1}^{\alpha}.
\end{equation*}
It follows from (\ref{I=J}) that 
\begin{equation*}
\frac{J_{n}^{\alpha}}{(\sharp V(G_n))^2} \gtrsim \left(\frac{2^\alpha}{N}%
\right)^n \longrightarrow\infty \quad (n\rightarrow\infty),
\end{equation*}
which implies $\dim_{\mathrm{DH}}(\mathcal{G})\leq s.$

Let now $0<\alpha <s.$ We decompose the sum into internal and external pairs
and obtain by (\ref{d_n=2^-}) that 
\begin{equation*}
I_{n}^{\alpha }=\sum_{j\in \Sigma }\sum_{\substack{ u,v\in j\Sigma ^{n-1} 
\\ u\neq v}}\left( \widetilde{d}_{n}(u,v)\right) ^{-\alpha }+E_{n}^{\alpha
}=2^{\alpha }NI_{n-1}^{\alpha }+E_{n}^{\alpha },  
\end{equation*}%
where 
\begin{equation*}
E_{n}^{\alpha }:=\sum_{i\neq j}\sum_{u,v\in \Sigma ^{n-1}}\left( \frac{1}{2}%
h_{n-1}(u)+\frac{1}{2}h_{n-1}(v)+2^{-n}d_{1}(i,j)\right) ^{-\alpha }.
\end{equation*}

By Corollary 1, there exists $c>0$, independent of $r$ and $k$, such that
for all $1\leq r<k$, 
\begin{equation*}
\delta _{r,k}\geq c2^{r}.
\end{equation*}%
For $v\in \Sigma ^{k}\setminus \{\omega |_{k}\}$, let $m=|v\wedge (\omega|_k) |$,
where $0\leq m<k$. If $m<k-1$, then by (\ref{d_n(x|n}), 
\begin{equation*}
h_{k}(v)\geq 2^{-k}\delta _{k-m-1,k}\geq \frac{c}{2}2^{-m}.
\end{equation*}%
If $m=k-1$, then 
\begin{equation*}
h_{k}(v)=2^{-k}d_{1}(v_{k},\omega _{k})\geq \frac{d_{\min }}{2}2^{-m},
\end{equation*}%
where $d_{\min }:=\min_{i\neq j}d_{1}(i,j)>0.$ Therefore, 
\begin{equation}
h_{k}(v)\geq c_{0}2^{-\left\vert v\wedge (\omega |_{k})\right\vert}\text{ for all }v\in \Sigma
^{k}\setminus \{\omega |_{k}\},  \label{hk-v}
\end{equation}%
where $c_{0}=\min \{c/2,d_{\min }/2\}$ is independent of $k$.

Note that 
\begin{equation*}
\sharp \{v\in \Sigma ^{k}:\left\vert v\wedge (\omega |_{k})\right\vert
=m\}=(N-1)N^{k-m-1}\text{ }(0\leq m<k),
\end{equation*}%
by (\ref{hk-v}) and $2^{\alpha }/N<1$ one has 
\begin{equation}
\sum_{\substack{ u\in \Sigma ^{n-1}  \\ u\neq \omega |_{n-1}}}\left(
h_{n-1}(u)\right) ^{-\alpha }\lesssim N^{n-1}\sum_{m=0}^{n-2}\left( \frac{%
2^{\alpha }}{N}\right) ^{m}\lesssim N^{n-1}.  \label{h-sum}
\end{equation}

To avoid the singular term $h_{n-1}(\omega |_{n-1})^{-\alpha},$ we decompose $%
E_{n}^{\alpha }$ into the three disjoint cases%
\begin{equation*}
u\neq \omega |_{n-1},\qquad u=\omega |_{n-1},v\neq \omega |_{n-1},\qquad
u=v=\omega |_{n-1}.
\end{equation*}%
Write $w=\omega |_{n-1}.$ It follows from (\ref{h-sum}) that 
\begin{align*}
E_{n}^{\alpha }& =\sum_{i\neq j}\left( \sum_{\substack{ u\neq w \\ u,v\in
\Sigma ^{n-1}}}+\sum_{\substack{ u=w,v\neq w \\ v\in \Sigma ^{n-1}}}%
+\sum_{u,v=w}\right) \left( \frac{1}{2}h_{n-1}(u)+\frac{1}{2}%
h_{n-1}(v)+2^{-n}d_{1}(i,j)\right) ^{-\alpha } \\
& \leq N(N-1)N^{n-1}2^{\alpha }\sum_{\substack{ u\in \Sigma ^{n-1} \\ u\neq w
}}h_{n-1}(u)^{-\alpha } \\
& \quad +N(N-1)2^{\alpha }\sum_{\substack{ v\in \Sigma ^{n-1} \\ v\neq w}}%
h_{n-1}(v)^{-\alpha }+2^{\alpha n}\sum_{i\neq j}d_{1}(i,j)^{-\alpha } \\
& \lesssim N^{2n-2}.
\end{align*}

Thus $I_{n}^{\alpha }$ satisfies the recurrence 
\begin{equation*}
I_{n}^{\alpha }\leq 2^{\alpha }NI_{n-1}^{\alpha }+CN^{2n-2}
\end{equation*}%
for some constant $C>0$ independent of $n$. Dividing by $N^{2n}$ gives 
\begin{equation*}
\frac{I_{n}^{\alpha }}{N^{2n}}\leq \frac{2^{\alpha }}{N}\frac{%
I_{n-1}^{\alpha }}{N^{2(n-1)}}+\frac{C}{N^{2}}.
\end{equation*}%
Since $2^{\alpha }/N<1$, it follows that $I_{n}^{\alpha }\lesssim N^{2n}.$ 
Together with (\ref{I=J}), we obtain 
\begin{equation*}
\varliminf_{n\rightarrow \infty }\frac{J_{n}^{\alpha }}{(\sharp V(G_{n}))^{2}%
}<\infty ,
\end{equation*}%
which implies $\dim _{\mathrm{DH}}(\mathcal{G})\geq s.$
\end{proof}

\bigskip

\section{Generalized Wiener indices}

\subsection{Asymptotic behavior via GHP convergence}

First, we investigate the asymptotic behavior of the generalized Wiener
index associated with the metric-preserving substitution system. Let 
\begin{equation*}
K_{\omega }:=\int_{(x,y)\in \Sigma ^{\mathbb{N}}\times \Sigma ^{\mathbb{N}}}%
\widetilde{d}(x,y)d\mu (x)d\mu (y).
\end{equation*}%
By Theorem \ref{th-GHP-1}, if $\omega \in \Sigma ^{\mathbb{N}}$ satisfies $%
\omega _{k}\neq \omega _{1}$ $(k\geq 2),$ then the normalized spaces $%
\mathcal{X}_{n}^{\omega }$ converge to $\mathcal{X}_{\omega }$ in the
Gromov--Hausdorff--Prokhorov sense. Consequently, 
\begin{equation}
\lim_{n\rightarrow \infty }\frac{\mathcal{W}_{n}}{2^{n}N^{2n}}=\frac{%
K_{\omega }}{2}.
\label{lim-K}
\end{equation}

\begin{proposition}
\label{prop-K_omega}It holds that 
\begin{equation*}
K_{\omega }=\frac{2(N-1)^{2}}{N(2N-1)}C_{\omega }.
\end{equation*}
\end{proposition}

\begin{proof}
Recall that 
\begin{equation*}
\widetilde{d}(x,y)=2^{-\left\vert x\wedge y\right\vert -1}(h(\sigma
^{\left\vert x\wedge y\right\vert +1}x)+h(\sigma ^{\left\vert x\wedge
y\right\vert +1}y)),
\end{equation*}%
then %
\begin{equation*}
K_{\omega }=\sum_{m=0}^{\infty }2^{-m-1}\int_{\left\vert x\wedge
y\right\vert =m}\left( h(\sigma ^{m+1}x)+h(\sigma ^{m+1}y)\right) d(\mu
\times \mu ).
\end{equation*}%
Let $A_{m}=\{(x,y):\left\vert x\wedge y\right\vert =m\}.$ Since $\mu $ is
the uniform Bernoulli measure, we obtain $(\mu \times \mu )(A_{m})=\frac{N-1}{N^{m+1}}.$ Moreover, by shift invariance of $\sigma ,$ the
random variables $\sigma ^{m+1}x$ and $\sigma ^{m+1}y$ are independent and
have the same distribution as $\mu .$ Therefore,%
\begin{equation*}
\int_{A_{m}}h(\sigma ^{m+1}x)d(\mu \times \mu )=(\mu \times \mu )(A_{m})%
\mathbb{E}_{\mu }[h],
\end{equation*}%
and the same holds for the term involving $y.$ Hence 
\begin{equation*}
K_{\omega }=2\mathbb{E}_{\mu }[h]\sum_{m=0}^{\infty }2^{-m-1}\frac{N-1}{N^{m+1}}=\frac{2(N-1)}{2N-1}\mathbb{E}_{\mu }[h].
\end{equation*}%
Finally, by writing 
\begin{equation*}
h(z)=C_{\omega }\sum_{k=1}^{\infty }2^{-k}\mathbf{1}_{\{z_{k}\neq \omega
_{k}\}},
\end{equation*}%
we have 
\begin{equation*}
\mathbb{E}_{\mu }[h]=C_{\omega }\sum_{k=1}^{\infty }2^{-k}\mu (z_{k}\neq
\omega _{k})=C_{\omega }\frac{N-1}{N}.
\end{equation*}%
Substituting this into the previous expression implies%
\begin{equation*}
K_{\omega }=\frac{2(N-1)^{2}}{N(2N-1)}C_{\omega }.
\end{equation*}
\end{proof}

\begin{proof}[Proof of Theorem \protect\ref{th-W}]
From (\ref{lim-K}), it suffices to determine the integral $\mathcal{D}_{\omega
}.$ Since $\mu _{\omega }=\pi _{\sharp }\mu ,$ we have 
\begin{equation*}
(\pi \times \pi )_{\sharp }(\mu \times \mu )=\mu _{\omega }\times \mu
_{\omega }\text{ and }d_{\omega }\circ (\pi \times \pi )=\widetilde{d}.
\end{equation*}%
Hence, by the change of variable formula for push-forward measures and
Proposition \ref{prop-K_omega}, 
\begin{equation*}
\mathcal{D}_{\omega }=K_{\omega }=\frac{2(N-1)^{2}}{N(2N-1)}C_{\omega }.
\end{equation*}
\end{proof}

\subsection{Closed form}

\label{subsection-closed}

The GHP convergence determines the leading asymptotic of $\mathcal{W}_{n}.$
To obtain a precise formula, it remains to analyse the recursive
contributions coming from the substitution structure, which are not captured
by the limiting metric space.

\begin{theorem}
\label{th-W_n}Let $\{(G_{n},d_{n})\}_{n\geq 1}$ be a metric-preserving $%
\omega $-substitution graph system. Then for each $n\geq 1,$ 
\begin{eqnarray*}
\mathcal{W}_{n} &=&\frac{2(N-1)^{2}}{N(2N-1)}2^{n}N^{2n}%
\sum_{k=2}^{n-2}2^{-k}D_{k}+(N-1)N^{n-1}\sum_{k=1}^{n-1}(n-k)N^{n-k}%
\varepsilon _{k} \\
&&+N^{n-1}\frac{N^{n}-1}{N-1}\mathcal{W}_{1}+(N-1)N^{n}%
\sum_{k=2}^{n-1}(n-k)(N^{n-1}-N^{n-k})D_{k} \\
&&-N^{2n-1}\sum_{k=2}^{n-2}\left( (n-k+1)N-(n-k+2)\right) D_{k}-\frac{1}{2N-1%
}N^{n}\sum_{k=2}^{n-2}N^{k}D_{k},
\end{eqnarray*}%
where $\varepsilon _{k}=\sum_{v\in \Sigma }d_{1}(v,\omega _{k}).$
\end{theorem}

\begin{corollary}
By specializing the distance function in Theorem \ref{th-W_n} to the geodesic distance and the effective resistance distance, one immediately obtains the
corresponding closed-form formulas for the Wiener index and the Kirchhoff
index, respectively.
\end{corollary}

Let 
\begin{equation*}
\mathcal{H}_{n}:=\sum_{u\in \Sigma ^{n}}d_{n}(u,\omega |_{n})
\end{equation*}%
be the sum of distances from all vertices to the interface vertex $\omega
|_{n}.$

Using Lemma \ref{le-delta}, we derive a closed form for the internal
distance sum $\mathcal{H}_{n}.$

\begin{proposition}
\label{le-S}It holds that for each $n>1,$ 
\begin{eqnarray*}
\mathcal{H}_{n} &=&N^{n}\sum_{k=1}^{n}N^{-k}\varepsilon
_{k}+(N-1)N^{n-1}\sum_{k=1}^{n-1}(n-k)N^{-k}\varepsilon _{k} \\
&&+\frac{N-1}{N}N^{n}2^{n+1}\sum_{k=2}^{n-1}2^{-k}D_{k}-(N-1)N^{n-1}%
\sum_{k=2}^{n-1}(n-k+2)D_{k} \\
&&+N^{n}\sum_{k=2}^{n}(1-N^{-k+1})D_{k}+(N-1)N^{n}%
\sum_{k=2}^{n-1}(n-k)(N^{-1}-N^{-k})D_{k}.
\end{eqnarray*}
\end{proposition}

\begin{proof}
For $1\leq k\leq n,$ let 
\begin{equation*}
\mathcal{H}_{k,n}:=\sum_{u\in \Sigma ^{k}}d_{k}(u,\sigma ^{n-k}(\omega
|_{n})).
\end{equation*}%
Then $\mathcal{H}_{k,k}=\mathcal{H}_{k}$ for all $k\geq 1.$ Using Claim \ref%
{clm-3}, we obtain that 
\begin{eqnarray}
\mathcal{H}_{k,n} &=&\sum_{j\neq \omega _{n-k+1}}\sum_{u\in j\Sigma
^{k-1}}d_{k}(u,\sigma ^{n-k}(\omega |_{n}))+\sum_{u\in \omega _{n-k+1}\Sigma
^{k-1}}d_{k}(u,\sigma ^{n-k}(\omega |_{n}))  \notag \\
&=&(N-1)\sum_{u\in \Sigma ^{k-1}}\left( d_{k-1}(u,\omega
|_{k-1})+d_{k-1}(\omega |_{k-1},\sigma ^{n-k+1}(\omega |_{n}))\right)  \notag
\\
&&+\sum_{j\neq \omega _{n-k+1}}\sum_{u\in \Sigma ^{k-1}}d_{1}(j,\omega
_{n-k+1})+\sum_{v\in \Sigma ^{k-1}}d_{k-1}(v,\sigma ^{n-k+1}(\omega |_{n})) 
\notag \\
&=&(N-1)\mathcal{H}_{k-1}+(N-1)N^{k-1}\delta _{k-1,n}+N^{k-1}\varepsilon
_{n-k+1}+\mathcal{H}_{k-1,n}.  \label{H_k,n}
\end{eqnarray}%
With $\mathcal{H}_{1,n}=\sum_{u\in \Sigma _{1}}d_{1}(u,\sigma ^{n-1}(\omega|
_{n}))=\varepsilon _{n},$ recursive equation (\ref{H_k,n}) gives us%
\begin{equation*}
\mathcal{H}_{n}=(N-1)\sum_{k=1}^{n-1}\mathcal{H}_{k}+\sum_{k=1}^{n}N^{n-k}%
\varepsilon _{k}+(N-1)\sum_{k=1}^{n-1}N^{k}\delta _{k,n}.
\end{equation*}%
For $k\geq 2,$ let $\Psi (k)=\mathcal{H}_{k}-N\mathcal{H}_{k-1}.$ Then 
\begin{equation*}
\mathcal{H}_{n}=N^{n-1}\mathcal{H}_{1}+\sum_{k=2}^{n}N^{n-k}\Psi (k)=N^{n-1}%
\mathcal{\varepsilon }_{1}+\sum_{k=2}^{n}N^{n-k}\Psi (k).
\end{equation*}%
Substituting (\ref{delta_k,n-new}), the explicit expression of $\delta
_{k,n},$ into $\Psi (k)$, we obtain the explicit formula for $\mathcal{H}%
_{n}.$
\end{proof}

\begin{remark}
With the notation introduced in Section 4, we have $\frac{\mathcal{H}_{n}}{%
2^{n}N^{n}}\rightarrow \mathbb{E}_{\mu }[h]$ as $n\rightarrow \infty .$
\end{remark}

Combining the result of internal distance sum $\mathcal{H}_{n}$ with the decomposition
of $\mathcal{W}_{n},$ we now return to Theorem \ref{th-W_n}.

\begin{proof}[Proof of Theorem \protect\ref{th-W_n}]
It follows from Claim \ref{clm-3} and the self-similarity of $G_{n}$ that%
\begin{eqnarray*}
\mathcal{W}_{n} &=&\sum_{j=1}^{N}\sum_{\{u,v\}\in j\Sigma
^{n-1}}d_{n}(u,v)+\sum_{1\leq i<j\leq N}\sum_{\substack{ u\in i\Sigma ^{n-1}
\\ v\in j\Sigma ^{n-1}}}d_{n}(u,v) \\
&=&\sum_{j=1}^{N}\sum_{\{u,v\}\in \Sigma ^{n-1}}d_{n-1}(u,v) \\
&&+\sum_{1\leq i<j\leq N}\sum_{\substack{ u\in i\Sigma ^{n-1} \\ v\in
j\Sigma ^{n-1}}}\left( d_{n-1}(\sigma u,\omega
|_{n-1})+d_{1}(u_{1},v_{1})+d_{n-1}(\omega |_{n-1},\sigma v)\right)  \\
&=&N\mathcal{W}_{n-1}+\sum_{1\leq i<j\leq N}\sum_{\substack{ u\in i\Sigma
^{n-1} \\ v\in j\Sigma ^{n-1}}}d_{1}(u_{1},v_{1}) \\
&&+\sum_{1\leq i<j\leq N}\sum_{\substack{ u\in \Sigma ^{n-1} \\ v\in j\Sigma
^{n-1}}}d_{n-1}(u,\omega |_{n-1})+\sum_{1\leq i<j\leq N}\sum_{\substack{ %
u\in i\Sigma ^{n-1} \\ v\in \Sigma ^{n-1}}}d_{n-1}(\omega |_{n-1},v) \\
&=&N\mathcal{W}_{n-1}+N^{2(n-1)}\mathcal{W}_{1}+(N-1)N^{n}\mathcal{H}_{n-1}.
\end{eqnarray*}%
Substituting the expression of $\mathcal{H}_{n-1}$ into the above recursion of 
$\mathcal{W}_{n}$ and iterating this recursion, we obtain the explicit
expression for $\mathcal{W}_{n}$.
\end{proof}

\bigskip

\section{Conclusions}

In this paper, we introduced a new class of substitution systems, called $%
\omega $-substitution graph systems, in which the interface vertices are
determined by an infinite sequence $\omega $. This construction provides a framework for generating non-homogeneous self-similar graph
sequences. Under suitable assumptions on $\omega $, we established the
Gromov--Hausdorff--Prokhorov convergence of the associated metric measure
spaces, proved the consistency between the discrete and continuous fractal
dimensions, and derived asymptotic as well as closed-form formulas for a
class of generalized Wiener indices.

A central feature of the present model is the consistency between its
discrete and continuous dimensional structures. More precisely,
\begin{equation*}
  \dim_{\mathrm{DH}}(\mathcal G)
=
\dim_{\mathrm{H}}(X_\omega)
=
\dim_{\mathrm{B}}(X_\omega)
=
\frac{\log N}{\log 2}.
\end{equation*}
Thus, for the $\omega$-substitution graph systems considered here, the
fractal dimension of the scaling limit can be recovered directly from the
approximating graph sequence. Together with the GHP convergence and the
asymptotic formulas for distance-based indices, this provides a unified
discrete-to-continuum description of the geometry of these systems.

Several directions remain open for future research, include investigating random $\omega $-substitution graph systems
and their probabilistic scaling limits, as well as spectral properties of
the associated graph sequences and limiting spaces. The present framework provides a foundation for further studies of
substitution graph systems from the perspective of metric geometry and
fractal analysis.

\bigskip

%
\section*{Conflict of interest}
The authors declare that they have no financial or non financial interests relevant to the content of this article

\section*{Funding}
The research is partially supported by the Natural Science Foundation of China (No. 12301107) and the Natural Science 
Foundation of Shandong Province, China (No. ZR202209010046). 


\end{document}